\documentclass[10pt]{amsart}
\newtheorem{theorem}{Theorem}[section]
\newtheorem{lemma}[theorem]{Lemma}
\newtheorem{proposition}[theorem]{Proposition}
 
\theoremstyle{definition}

\newtheorem{remark}[theorem]{Remark}

\newcommand{\id}{\operatorname{id}} 
 
\newcommand{\Ker}{\text{Ker\,}}

\newcommand{\Ad}{\operatorname{Ad}}

\newcommand{\Aut}{\text{Aut}}

\newcommand{\nc}{\newcommand}
\nc{\Symm}{{\on{Sym}}}

\newcommand{\on}{\operatorname}   
\newcommand{\eps}{\varepsilon}
 \nc{\cE}{{\cal E}}

\newcommand{\cG}{\mathcal{G}}

\nc{\SL}{{\mathfrak sl}}
\nc{\HH}{{\mathfrak h}}
\newcommand{\g}{{\mathfrak{g}}}

\renewcommand{\k}{{\mathfrak{k}}}

\newcommand{\SG}{{\mathfrak{S}}}

\nc{\wh}{\widehat}\nc{\wt}{\widetilde}

\newcommand{\la}{{\lambda}}

\newcommand{\ben}{\begin{enumerate}}
\newcommand{\een}{\end{enumerate}}

\newcommand{\cO}{{\mathcal O}}

\newcommand{\cR}{{\mathcal R}}

\newcommand{\kk}{{\bf k}}

\newcommand{\QQ}{{\mathbb{Q}}}

\newcommand{\RR}{{\bold R}}
\newcommand{\LL}{{\bold L}}

\renewcommand{\Ad}{{\operatorname{Ad}}}

\begin{document}

\title[Comparison of Poisson structures and Poisson-Lie dynamical 
$r$-matrices] {Comparison of Poisson structures and Poisson-Lie 
dynamical $r$-matrices}

\begin{abstract} 
We construct a Poisson isomorphism between the formal Poisson 
manifolds $\g^*$ and $G^*$, where $\g$ is a finite dimensional 
quasitriangular Lie bialgebra. Here $\g^*$ is equipped with its 
Lie-Poisson (or Kostant-Kirillov-Souriau) structure, and $G^*$ with 
its Poisson-Lie structure. We also quantize the Poisson-Lie dynamical 
$r$-matrices of Balog-Feh\'er-Palla. 
\end{abstract}


\author{Benjamin Enriquez}
\address{IRMA (CNRS), 7 rue Ren\'e Descartes, F-67084 Strasbourg, France}
\email{enriquez@@math.unistra.fr}

\author{Pavel Etingof}
\address{Department of Mathematics, Massachusetts Institute of Technology,
Cambridge, MA 02139, USA}
\email{etingof@@math.mit.edu}

\author{Ian Marshall}
\address{Department of Mathematics, National Research University Higher School of Economics, Usacheva str. 6, 119048, Moscow, Russia}
\email{imarshall@@hse.ru}

\maketitle

\section*{Introduction and main results}

We construct Poisson isomorphisms between the formal Poisson 
manifolds $\g^*$ and $G^*$, where $\g$ is a finite dimensional 
quasitriangular Lie bialgebra (Thm. \ref{thm:main}). Here $\g^*$ is equipped with its 
Lie-Poisson (or Kostant-Kirillov-Souriau) structure, and $G^*$ with 
its Poisson-Lie structure. 

Thm. \ref{thm:main} may be viewed as a generalization of the 
formal version of \cite{GW} (later reproved in \cite{Al}, 
Thm. 1, and \cite{Bo}), 
where Ginzburg and Weinstein construct a Poisson 
diffeomorphism between the real Poisson manifolds $\k^*$ and 
$K^*$, where $K$ is a compact Lie group and $\k$ is its Lie 
algebra. It can also be viewed as a new result in the subject of
linearization of Poisson structures; e.g., in contrast with our result, 
it has been shown in \cite{Ch} that not all Poisson structures 
on Poisson-Lie groups are linearizable. 

We give two proofs of Thm. \ref{thm:main}. The first proof relies on a nondegeneracy assumption 
and is geometric. Namely, it relies on the construction of a map $g(\lambda) :  \g^* \to G$
satisfying a differential equation, which is achieved by using the theory of the classical 
Yang-Baxter equation and gauge transformations. Note that a geometric proof of the same result, not relying on the nondegeneracy assumption, is given in \cite{AM2}.
The second proof relies of the
theory of quantization of Lie bialgebras. 

We then apply Thm. \ref{thm:main} to the quantization of Poisson-Lie 
dynamical $r$-matrices introduced in \cite{FM1} (see Thm. \ref{thm:2}). 

We now describe our results in more detail. 

\subsection{Comparison of Poisson structures: statement and proofs of Thm. \ref{thm:main}}

\subsubsection{Formulation of Thm. \ref{thm:main}}\label{sect:announcement}

Let $(\g,r)$ be a finite dimensional quasitriangular Lie bialgebra
over a field $\kk$ of characteristic $0$. Recall that this means that 
$\g$ is a Lie algebra, $r\in \g^{\otimes 2}$, 
$t:= r+ r^{2,1}$ is invariant, and $\on{CYB}(r) := 
[r^{1,2},r^{1,3}] + [r^{1,2},r^{2,3}] + [r^{1,3},r^{2,3}] = 0$. The Lie
cobracket on $\g$ is defined by $\delta(x) := [x\otimes 1+1\otimes x,r]$
for any $x\in \g$. 

The Lie bracket on $\g$ defines a linear Poisson structure on $\g^*$; we also 
denote by $\g^*$ the formal neighborhood of the origin in this vector space,
which is then a formal Poisson manifold.  
On the other hand, we denote by $G^*$ the formal Poisson-Lie group with Lie 
bialgebra $\g^*$. Our main result is: 

\begin{theorem} \label{thm:main}
The formal Poisson manifolds $\g^*$ and $G^*$ are isomorphic.  
\end{theorem}

\subsubsection{Structure of the proof of Thm. \ref{thm:main}}

Thm. \ref{thm:main} will be proved in two ways: (a) under the assumption that $t\in S^2(\mathfrak g)$
is nondegenerate as a consequence of Props. \ref{prop:1} and \ref{prop:2} (proved geometrically in 
Sect. \ref{sect:geom} and using quantization in Sects. \ref{subsect:toto} and \ref{subsect:titi}); 
(b) unconditionally in Sect. \ref{section1:2.5}.   

\subsubsection{The statements leading to the proof of Thm. \ref{thm:main}}\label{subsect:marseille}

We now formulate Props. \ref{prop:1} and \ref{prop:2}. 
Denote by $G$ the formal group with Lie algebra $\g$ and by 
$\on{Map}_0(\g^*,G)$ the space of formal maps $g : \g^* \to G$, 
such that $g(0)=1$; this is the space of maps of the form 
$e^{x(\lambda)}$, where\footnote{We denote by $\wh S(\g)$ the degree completion 
$\wh\oplus_{k\geq 0} S^k(\g)$ of the symmetric algebra $S(\g)$, 
and set $\wh S(\g)_{>i} = \wh\oplus_{k>i} S^k(\g)$.}   $x(\lambda)\in \g\otimes \wh S(\g)_{>0}$.
The set 
$\on{Map}_0(\g^*,G)$ has a group structure, defined by 
$(g_1 * g_2)(\lambda) := g_2(\on{Ad}^*(g_1(\lambda))(\lambda)) 
g_1(\lambda)$. Its subspace of all maps $g(\lambda)$ such that 
\begin{equation} \label{eq:map:ham}
g_1^{-1}d_2(g_1)(\lambda) - g_2^{-1}d_1(g_2)(\lambda) + 
\langle \id\otimes \id\otimes \lambda , 
[g_1^{-1} d_3(g_1)(\lambda), g_2^{-1} d_3(g_2)(\lambda)] \rangle  = 0
\end{equation} 
is a subgroup $\on{Map}_0^{ham}(\g^*,G)$. 
Here $g_1^{-1}d_2(g_1)(\lambda) = 
\sum_i g^{-1}\partial_{\eps^i}g(\lambda) \otimes e_i$ is viewed as a 
formal function $\g^* \to \g^{\otimes 2}$,  $(\eps^i)$, $(e_i)$ are dual bases
of $\g^*$ and $\g$, $g_i^{-1}d_j(g_i) = (g_1^{-1}d_2(g_1))^{i,j}$ 
and $\partial_\xi g(\lambda) = (d/d\eps)_{|\eps=0} g(\lambda + \eps\xi)$.  
We will also denote by $g_{13}^{-1}d_2(g_{13})$ the same quantity, viewed 
as an element of $\g^{\otimes 2} \otimes \wh S(\g)$. 

We have a group morphism 
$\theta : \on{Map}_0^{ham}(\g^*,G) \to \on{Aut}_1(\g)^{op}$ 
to the group (with opposite structure) of Poisson automorphisms 
of $\g^*$ with differential at $0$ equal to the identity, 
taking $g(\lambda)$ to the automorphism 
$\lambda\mapsto \Ad^*(g(\lambda))(\lambda)$. 

\begin{proposition}  \label{prop:1}
There exists a formal map $g(\lambda) \in \on{Map}_0(\g^*,G)$, 
such that 
\begin{equation} \label{dg=rho}
(g_1)^{-1}d_2(g_1) - (g_2)^{-1}d_1(g_2) + \Ad(g\otimes g)^{-1}(r_0) 
 + \langle \id\otimes \id \otimes \lambda, 
[(g_1)^{-1}d_3(g_1),(g_2)^{-1}d_3(g_2)]\rangle 
= \rho_{\on{AM}}. 
\end{equation}
(identity of 
formal\footnote{We view $\wedge^n(\g)$ as a subspace of $\g^{\otimes n}$.}  
maps $\g^*\to\wedge^2(\g)$). 
Here $r_0 = (r-r^{2,1}) / 2$, and $\rho_{\on{AM}}$ is the 
Alekseev-Meinrenken $r$-matrix (\cite{AM,BDF}) given by 
$$
\rho_{\on{AM}}(\lambda) = (\id\otimes \varphi(\on{ad} \lambda^\vee))(t), 
$$ 
where $\lambda^\vee = (\lambda\otimes \id)(t)$ and 
$\varphi(z) := -{1\over z} + {1\over 2} \on{cotanh}{z\over 2}$.  

The group $\on{Map}_0^{ham}(\g^*,G)$ acts simply and transitively 
on the space of solutions $g(\lambda)$ of (\ref{dg=rho}), as follows: 
$(\alpha * g)(\lambda) = 
g(\Ad^*(\alpha(\lambda))(\la)) \alpha(\la)$.
\end{proposition}

\begin{proposition} \label{prop:2}
Assume that $t\in S^2(\mathfrak g)$ is nondegenerate. 
Let $g(\lambda)\in \on{Map}_0(\g^*,G)$ be as in Prop. \ref{prop:1}. 
There exists a unique isomorphism $g^*(\lambda) : \g^* \to G^*$, defined by 
the identity 
$$
g(\lambda) e^{\lambda^\vee} g(\lambda)^{-1}
= L(g^*(\lambda))R(g^*(\lambda))^{-1}. 
$$
Here $L,R : G^* \to G$ are the formal group morphisms corresponding 
to the Lie algebra morphisms $L,R: \g^* \to \g$  given by 
$L(\lambda) := (\lambda\otimes \id)(r)$, 
$R(\lambda) := -(\lambda\otimes \id)(r^{2,1})$.  

In other words, we have (non-Poisson) formal manifold isomorphisms 
$\g^* \stackrel{a}{\to} G \stackrel{b}{\gets} G^*$, 
$a(\lambda) = g(\lambda) e^{\lambda^\vee} g(\lambda)^{-1}$, 
$b(g^*) = L(g^*)R(g^*)^{-1}$, and $g^*(\lambda) 
= b^{-1} \circ a(\lambda)$. 

The isomorphism $(\alpha * g)^*(\la) : \g^* \to G^*$
corresponding to $(\alpha * g)(\la)$ is such that 
$(\alpha * g)^*(\la)  = g^*(\theta(\alpha)(\la))$. 
\end{proposition}

Props. \ref{prop:1} and \ref{prop:2} immediately imply Thm. \ref{thm:main} 
under the assumption that $t\in S^2(\mathfrak g)$ is nondegenerate. 

These propositions also imply that 
the set of all isomorphisms $\g^*\to G^*$ that they allow to construct 
is a 
principal homogeneous space under the image of 
$\theta : \on{Map}_0^{ham}(\g^*,G)
\to \Aut_1(\g^*)$. When $\g$ is semisimple, any derivation 
$\g\to S^k(\g)$ is inner, so $\theta$ is surjective. So in that 
case, our construction yields all 
the Poisson isomorphisms $\g^*\to G^*$ taking $0$ to $1$ and 
with differential at $0$ equal to the identity. 

\subsection{Quantization of Poisson-Lie dynamical $r$-matrices}

In \cite{BFP,FM1}, there was introduced the Poisson-Lie dynamical 
$r$-matrix 
$$
\rho_{\on{FM}}(g^*) = \Big( \id \otimes \big( 
\nu {{\id + a(g^*)^{2\nu}}\over{\id- a(g^*)^{2\nu}}} 
- {1\over 2}{{\id+a(g^*)}\over{\id-a(g^*)}}\big) \Big)(t). 
$$
Here $a(g^*) : G^* \to \on{GL}(\g)$ is defined by $a(g^*)
= \Ad( L(g^*)R(g^*)^{-1})$ and $\nu$ is a fixed element in $\kk$. 
Then $\rho_{\on{FM}}$ is a formal map $G^* \to \wedge^2(\g)$, 
i.e., an element of $\wedge^2(\g) \otimes \cO_{G^*}$. 
It is a Poisson-Lie classical dynamical $r$-matrix for 
$(G^*,\g,Z_\nu)$, where $Z_{\nu} = (\nu^2 - {1\over 4})[t^{1,2},t^{2,3}]$
(see also \cite{EEM}).  

A quantization of this dynamical $r$-matrix is the data of 
a quantization $U_\hbar(\g)$ of the Lie bialgebra $\g$, 
and a pair $(\bar J,\bar \Phi)$, where: 

$\bullet$  $\bar \Phi$ is an associator for $U_\hbar(\g)$, i.e., 
$\bar \Phi \in U_\hbar(\g)^{\wh\otimes 3}$ satisfies
\footnote{We denote by $\wh\otimes$ the topological tensor product, defined 
as follows: if $V,W$ are topological vector spaces of the form 
$V = V_0[[x_1,\ldots,x_n]]$, $W = W_0[[y_1,\ldots,y_m]]$, then 
$V\wh\otimes W := V_0\otimes W_0[[x_1,\ldots,y_m]]$.} 
\footnote{We set $\bar\Phi^{12,3,4} = 
(\Delta \otimes \id\otimes \id)(\bar\Phi)$, etc.}
the pentagon equation $\bar \Phi^{2,3,4} \bar\Phi^{1,23,4} \bar\Phi^{1,2,3} 
= \bar\Phi^{1,2,34} \bar\Phi^{12,3,4}$, 
the invariance condition $[\bar\Phi,(\Delta\otimes \id)\circ \Delta(a)] = 0$
for any $a\in U_\hbar(\g)$, 
and the expansion $\bar\Phi = 1 + \hbar^2 \phi_2$, where $\phi_2\in
U_\hbar(\g)^{\wh\otimes 3}$ satisfies 
$\sum_{\sigma\in \SG_3} \on{sign}(\sigma) \sigma((\phi_2)_{|\hbar=0}) 
= Z_\nu$; 

$\bullet$ $\bar J\in
U_\hbar(\g)^{\wh\otimes 2} \wh\otimes U_\hbar(\g)'$ satisfies the 
dynamical twist equation $\bar J^{2,3,4}\bar J^{1,23,4} \bar\Phi^{1,2,3}
= \bar J^{1,2,34} \bar J^{12,3,4}$, 
the invariance condition $[\bar J,(\Delta\otimes \id)\circ \Delta(a)] = 0$
for any $a\in U_\hbar(\g)$, and the expansion
$\bar J = 1 + \hbar \bar j_1$, where $\bar j_1\in U_\hbar(\g)^{\wh\otimes 2} 
\wh\otimes U_\hbar(\g)'$ satisfies 
$(\bar j_1 - \bar j_1^{2,1,3})_{|\hbar=0} = \rho_{\on{FM}}$. 

\begin{theorem} \label{thm:2}
Such a quantization exists (see Thm. \ref{thm:end}). 
\end{theorem}

\section{Construction of isomorphisms $\g^* \to G^*$}\label{section:1}

We give two families of proofs of the statements announced in Sect. \ref{sect:announcement}. The first family 
is geometric and is contained in Sect. \ref{sect:geom}. The second family relies on quantization of Lie bialgebras and is 
contained in Sect. \ref{sect:quant}. 

Recall that the statements from Sect. \ref{sect:announcement} comprise Thm. \ref{thm:main}, 
Prop. \ref{prop:1} and Prop. \ref{prop:2}. 
It is noted in Sect. \ref{subsect:marseille} that Thm. \ref{thm:main} follows from Props. \ref{prop:1} and 
\ref{prop:2} if $t$ is nondegenerate. We give two proofs of these propositions: (a) a Poisson geometric 
proof (Sect. \ref{sect:geom}), and (b) a proof based on the theory of quantization of Lie bialgebras (Sects. 
\ref{subsect:toto} and \ref{subsect:titi}).  In Sect. \ref{section1:2.5}, we give an unconditional proof of 
Thm. \ref{thm:main}. 
 
\subsection{Geometric construction} \label{sect:geom}

\subsubsection{Construction of $g(\lambda)$ (proof of Prop. \ref{prop:1})}

One checks that $\on{Map}^{ham}_0(\g^*,G)$ is a prounipotent Lie
group with Lie algebra $\{\alpha\in \g\otimes \wh S(\g)_{\geq 1} 
| \on{Alt} \circ d(\alpha) = 0\}$. This Lie algebra is isomorphic to 
$(\wh S(\g)_{>1},\{-,-\})$ under $d : f\mapsto d(f)$. 
\footnote{For $\rho = \alpha \otimes f\in\wedge^n(\g)\otimes\wh S(\g)$, 
we set $d(\rho) := \sum_i \alpha\otimes e_i \otimes 
(d/d\eps)_{|\eps=0} f(\lambda + \eps \eps^i)$ and if 
$\xi\in \wedge^{n-1}(\g) \otimes \g$,   
$\on{Alt}(\xi\otimes f) := (\xi + \xi^{2,\ldots,n,1} + \cdots  
+ \xi^{n,1,\ldots,n-1})\otimes f$.} 

Let us denote by $\cG$ the set of all $g\in \on{Map}_0(\g^*,G)$
satisfying (\ref{dg=rho}). This is a subvariety of the proalgebraic 
variety $\on{Map}_0(\g^*,G)$. One checks that 
$(\wh S(\g)_{>1},\{-,-\})$ acts by vector fields on 
$\cG$, by 
\begin{equation} \label{action}
g^{-1}\delta_f(g) = \langle \id\otimes \id\otimes \lambda,
[d_3(f_2),g_{12}^{-1} d_3(g_{12})] \rangle - d_1(f_2) \in \g \otimes 
\wh S(\g)_{>0},  
\end{equation}
and that the right infinitesimal action of $\on{Map}_0(\g^*,G)$
on itself is given by the same formula. It follows that the right action of 
$\on{Map}_0(\g^*,G)$ on itself restricts to an action of 
$\on{Map}_0^{ham}(\g^*,G)$ on $\cG$. 

We now prove that if $\cG$ is nonempty, then $\on{Map}_0^{ham}(\g^*,G)$
acts simply and transitively on $\cG$. 

Let us show that the action is simple. If $g,g'\in\cG$ and $\alpha\in 
\on{Map}_0^{ham}(\g^*,G)$ are such that $g * \alpha = g'$, then let 
$a := \log(\alpha)$. Assume that $a\neq 0$ and let $n$ be the smallest integer
such that the component $a_n$ of $a$ in $\g\otimes S^n(\g)$ is $\neq 0$. Then 
$\log(g') - \log(g) = a_n$ modulo $\g\otimes \wh S(\g)_{>n}$, which implies 
that $g\neq g'$. 

Let us now prove that the action is transitive. Let $g,g'\in \cG$. Set $A:= 
\on{log}(g)$, $A' := \log(g')$. Then $A,A'\in \g\otimes 
\wh S(\g)_{>0}$. Assume that $A\neq A'$ and let $n$ be the 
smallest integer such that the component of $A'-A$ in $\g\otimes S^n(\g)$
is nonzero; we denote by $(A'-A)_n$ this component. Comparing 
equations (\ref{dg=rho}) for $g$ and $g'$, we get $\on{Alt} 
\circ d((A'-A)_n) = 0$. 
It follows that there exists $a\in S^{n+1}(\g)$, such that $(A'-A)_n 
= \on{Alt} \circ d(a)$, i.e., $(A'-A)_n \in \on{Lie}\on{Map}_0^{ham}(\g^*,G)$. 
Let $\exp_*$ be the exponential map of $\on{Map}_0^{ham}(\g^*,G)$
and $\alpha := \exp_*((A'-A)_n) \in \on{Map}_0^{ham}(\g^*,G)$. Then  
$\log(\alpha * g) = A(\Ad^*(\alpha(\lambda))(\lambda)) + (A'-A)_n$ 
modulo $\g\otimes \wh S(\g)_{\geq n+1}$; therefore the difference of 
logarithms
\footnote{Here the logarithm is the inverse of the ordinary exponential map 
$\g\otimes \wh S(\g) \to \on{Map}_0^{ham}(\g^*,G)$.} 
of $\alpha * g$ and $g'$ coincide modulo $\g\otimes 
\wh S(\g)_{\geq n+1}$. Working by successive
approximations, we construct $\beta\in \on{Map}_0^{ham}(\g^*,G)$ such that 
$\beta * g = g'$. This proves that the action is transitive. 

Let us now prove that $\cG$ is nonempty. Recall that $r_0 = (r-r^{2,1})/2$. 
If $g\in \on{Map}_0^{ham}(\g^*,G)$, set $(r_0)^g :=$ l.h.s. of 
(\ref{dg=rho}). 

\begin{lemma} \label{lemma:1}
Assume that $\log(g) = -r/2$ modulo $\g\otimes\wh S(\g)_{>1}$. 
Set $\rho := (r_0)^g$ and assume that $\rho = \rho_{inv} + \alpha$, 
where $\rho_{inv},\alpha\in
\wedge^2(\g) \otimes \wh S(\g)$, $\rho_{inv}$ is $\g$-invariant and 
$\alpha \in \g\otimes \wh S(\g)_{\geq n}$. Then 
$$
\on{CYB}(\rho) - \on{Alt}(d\rho) = Z^{1,2,3} 
- {1\over 2}[r^{1,234},\alpha^{2,3,4}]
$$
modulo $\g^{\otimes 3} \otimes \wh S(\g)_{\geq n+1}$.
\end{lemma}

{\em Proof of Lemma.} This statement can be proved directly. It can also be
viewed as the classical limit of the following statement. Let 
$U = U(\g)[[\hbar]]$, let $J\in (U^{\wh\otimes 2})^\times$, 
and let $\Phi := (J^{1,2}J^{12,3})^{-1} 
J^{2,3} J^{1,23}$. 
If $K\in (U^{\wh\otimes 2})^\times$ and we set 
$\bar\Phi := (K^{1,23})^{-1} (K^{2,3})^{-1} J^{1,2} K^{12,3}$, then 
$\bar\Phi$ satisfies 
$$
\bar\Phi^{2,3,4} \bar\Phi^{1,23,4} \Phi^{1,2,3} = 
(\bar\Phi^{2,3,4},(K^{1,234})^{-1}) \bar\Phi^{1,2,34} 
\bar\Phi^{12,3,4} (K^{123,4},\Phi^{1,2,3}). 
$$
Here $(a,b) = aba^{-1}b^{-1}$. The statement of the lemma is 
recovered when $J,K\in (U^{\wh\otimes 2})^\times$  
have the form $1-\hbar r/2 + o(\hbar)$, and $K$ is admissible 
with classical limit $g(\la)$: the contribution of 
$(\bar\Phi^{2,3,4},(K^{1,234})^{-1})$ is 
$- {1\over 2}[r^{1,234},\alpha^{2,3,4}]$, while 
the commutator $(K^{123,4},\Phi^{1,2,3})$ does not contribute 
(as the classical limit of $\Phi$ is proportional to $Z$ and 
hence invariant). \hfill \qed \medskip  

Let us now prove that $\cG$ is nonempty. We will construct 
a sequence $g_n\in \on{Map}_0(\g^*,G)$, such that 
$\rho_n := (r_0)^{g_n}$ satisfies $\rho_n = \rho_{\on{AM}}$
modulo $\wedge^3(\g) \otimes \wh S(\g)_{\geq n+1}$. 

If $n=0$, we set $g_0 := \exp(-r/2)$, then $\rho_0 = 0 = \rho_{\on{AM}}$
modulo $\wedge^3(\g) \otimes \wh S(\g)_{\geq 1}$. Assume that we determined 
$g_n$ and let us construct $g_{n+1}$. Set $\alpha := 
\rho_n - \rho_{\on{AM}}$, and let $\alpha_{n+1}$ be the 
component of $\alpha$ in $\wedge^2(\g) \otimes S^{n+1}(\g)$. 
Then Lemma \ref{lemma:1} implies that 
$$
\on{CYB}(\rho_{\on{AM}} + \alpha_{n+1}) - \on{Alt}(d(\rho_{\on{AM}} 
+ \alpha_{n+1})) = Z^{1,2,3} - {1\over 2}[r^{1,234},\alpha_{n+1}^{2,3,4}]
\on{\ mod\ }\wedge^3(\g) \otimes S^{n+1}(\g). 
$$
Since $\rho_{\on{AM}}$ satisfies the modified classical dynamical Yang-Baxter
equation, the component in $\wedge^3(\g) \otimes S^n(\g)$ of this identity
yields $\on{Alt}(d\alpha_{n+1}) = 0$, so that we have $\alpha_{n+1} =
\on{Alt}(d\beta)$ for some $\beta\in \g\otimes S^{n+2}(\g)$. Then we 
set $g_{n+1} := \exp_*(-\beta) * g_n$, $\rho_{n+1} := (r_0)^{g_{n+1}}$. 
Then $\rho_{n+1} - \rho_{\on{AM}} = \rho_n - \rho_{\on{AM}} 
- \on{Alt}(d\beta) = \alpha_{n+1} - \on{Alt}(d\beta) = 0$ modulo 
$\wedge^3(\g)\otimes \wh S(\g)_{\geq n+1}$. 
By successive approximations, we then construct $g$ such that 
$(r_0)^g = \rho_{\on{AM}}$. So $\cG$ is nonempty.   

\subsubsection{Poisson isomorphism $\g^* \to G^*$ (proof of Prop. 
\ref{prop:2})}

In this section, we show that the fact that $g(\lambda)$ 
satisfies (\ref{dg=rho}) implies that $\lambda \mapsto g^*(\lambda)$ 
is a Poisson isomorphism. We will freely use the formalism of differential
geometry, even though we work in the formal setup; the computations below make
sense (and prove the desired result) when working over an Artinian $\kk$-ring. 

According to \cite{STS}, $b : G^* \to G$ formal isomorphism if $t$ is nondegenerate and 
the image of the Poisson bracket $\{-,-\}_{G^*}$ on 
$G^*$ under $b$ is the Poisson bracket on $G$ 
$$
\{F,H\}_G(g) = 
\langle (d_{\RR} - d_\LL) F(g) \otimes d_\LL H(g), r \rangle
+ \langle (d_{\RR} - d_\LL) F(g) \otimes d_\RR H(g), r^{2,1} \rangle,  
$$
where $g\in G$, $F,H$ are functions on $G$, 
$d_\LL F(g), d_\RR F(g) \in \g^*$ are the left and right
differentials defined by 
$\langle d_\LL F(g) , a\rangle = (d/d\eps)_{|\eps=0} F(e^{\eps a}g)$, 
$\langle d_\RR F(g) , a\rangle = (d/d\eps)_{|\eps=0} F(ge^{\eps a})$ 
for any $a\in \g$. 

For $\xi\in \g^*$, define $F_\xi\in \cO_{G}$ by 
$F_\xi(g) = \langle \xi, \log(g) \rangle$.  
\begin{lemma}
$(d_\LL F_\xi)(g) = f({1\over 2} \on{ad}^*(\log g))(\xi)$ and 
$(d_\RR F_\xi)(g) = f(-{1\over 2} \on{ad}^*(\log g))(\xi)$, where 
$f(z) = z e^z / (\on{sinh} z)$, and $\on{ad}^*$ denotes the coadjoint action
of $\g$ on $\g^*$. 
\end{lemma}

{\em Proof.} Set $x := \log(g)$, let $a\in \g$, and set 
$\tilde a := (d/d\eps)_{|\eps=0} \log(e^{\eps a} e^x)$. 
The coefficient in $\eps$ of 
$e^{\eps a}e^x = e^{x + \eps\tilde a + O(\eps^2)}$ yields 
$ae^{x} = \sum_{n\geq 1} (n!)^{-1} \sum_{k=0}^{n-1} 
x^k \tilde a x^{n-1-k}$. Applying $\on{ad}(x)$ to this relation, we get 
$\tilde a = f(-{1\over 2} \on{ad} x)(a)$. 
Now $\langle d_\LL F_\xi(g), a\rangle = (d/d\eps)_{|\eps=0} F_\xi(e^{\eps a}g)
= \langle\xi,\tilde a \rangle = \langle \xi, f(-{1\over 2} \on{ad} x)(a)
\rangle = \langle f({1\over 2} \on{ad}^* x)(\xi), a\rangle$, which yields the
first formula. 
On the other hand, we have for any function $F$ on $G$, 
$d_\RR F(g) = \Ad^*(g^{-1})(d_\LL F(g))$, where $\on{Ad}^*$ is the coadjoint 
action of $G$ on $\g^*$, hence 
$d_\RR F_\xi(g) = \Ad^*(e^{-x}) f({1\over 2} \on{ad}^* x)(\xi)
= f(-{1\over 2} \on{ad}^* x)(\xi)$. 
\hfill \qed \medskip 

Then we get 
\begin{equation} \label{form:xi}
\{F_\xi,F_\eta\}_G(g) =  
\langle \on{ad}^*(x)(\xi) \otimes \on{ad}^*(x)(\eta) , 
(\id \otimes \varphi(\on{ad}(x))(t) + r_0 \rangle 
- \langle \on{ad}^*(x)(\xi) \otimes \eta, t\rangle, 
\end{equation}
where $x = \on{log}(g)$ and $\varphi$ and $r_0$ are as in Prop.  
\ref{prop:1}. 

On the other hand, let $\{-,-\}_G'$ be the image of the Poisson bracket 
$\{-,-\}_{\g^*}$ by $a : \g^* \to G$. Set $f_\xi(\lambda) := F_\xi\circ
a(\lambda)$, then $f_\xi(\lambda) = \langle \xi\otimes \lambda,
(\on{Ad}(g(\lambda))\otimes \id)(t) \rangle$. If $f\in \cO_{\g^*}$, 
$\la\in \g^*$, define $df(\la)\in\g$ by $\langle \alpha,df(\la)\rangle = 
(d/d\eps)_{|\eps=0} f(\la+\eps \alpha)$ for any $\alpha\in\g^*$. We have 
$\{f,g\}_{\g^*}(\la) = \langle \lambda, [df(\lambda),dg(\lambda)] \rangle$. 
Then $df_\xi(\la) = \langle \xi\otimes \id, A(\la)\rangle$, where 
$$
A(\la) = (\on{Ad}(g(\la)) \otimes \id)(t) 
+ [(d_2g_1)g_1^{-1}, \Ad(g(\la))(\la^\vee) \otimes 1].  
$$
So $$
\{f_\xi(\la),f_\eta(\la)\} = \langle \xi\otimes \eta\otimes \la,
[A^{1,3}(\la),A^{2,3}(\la)] \rangle. 
$$
This decomposes as the sum of four terms (we set $\bar x = 
\on{Ad}(g(\la))(\la^\vee)$): 

(a) $\langle \xi\otimes \eta\otimes \la, (\on{Ad}(g(\la))^{\otimes 2}
\otimes \id)([t^{13},t^{23}])\rangle = 
- \langle \xi\otimes \eta, (\on{ad}(\bar x) \otimes \id)(t)\rangle$; 

(b) $\langle \xi\otimes \eta\otimes \la, 
[[(d_3g_1)g_1^{-1},\Ad(g_1(\la))(\la^\vee)], \on{Ad}(g_2(\la))(t^{23})]\rangle 
= \langle \on{ad}^*(\bar x)^{\otimes 2}(\xi\otimes \eta), 
-\on{Ad}(g)^{\otimes 2}(g_1^{-1} d_2(g_1))\rangle$; 

(c) $\langle \xi\otimes\eta\otimes\la, 
[\on{Ad}(g_1(\la))(t^{13}), [(d_3g_2)g_2^{-1},
\on{Ad}(g_2(\la))(\la_2^\vee)]] \rangle
= \langle \on{ad}^*(\bar x)^{\otimes 2}(\xi\otimes \eta), 
\on{Ad}(g)^{\otimes 2}(g_2^{-1} d_1(g_2))\rangle$; 

(d) 
\begin{eqnarray*}
&\langle \xi\otimes\eta\otimes\la, \big[ 
[(d_3g_1)g_1^{-1},\on{Ad}(g_1(\la))(\la_1^\vee)], 
[(d_3g_2)g_2^{-1},\on{Ad}(g_2(\la))(\la_2^\vee)] \big] 
\\ & = \langle \on{ad}^*(\bar x)(\xi) \otimes 
\on{ad}^*(\bar x)(\eta) \otimes \la,
[(d_3g_1)g_1^{-1},(d_3g_2)g_2^{-1}] \rangle.
\end{eqnarray*}

We have $\{F_\xi,F_\eta\}'_G(g) = \{f_\xi,f_\eta\}_{\g^*}(a^{-1}(g))$; 
then $\bar x = x$, so $\{F_\xi,F_\eta\}'_G = \{F_\xi,F_\eta\}_G$
iff 
\begin{eqnarray*}
& \langle \on{ad}^*(x)(\xi) \otimes \on{ad}^*(x)(\eta), 
\on{Ad}(g)^{\otimes 2}(g_2^{-1} d_1(g_2) - g_1^{-1} d_2(g_1))
\\ & + \langle \id\otimes \id\otimes \la,[d_3(g_1)g_1^{-1}, 
d_3(g_2) g_2^{-1}]\rangle\rangle
- \big( \id\otimes \varphi(\on{ad}(x)))(t) - r_0 \rangle = 0, 
\end{eqnarray*}
for which a sufficient condition is that $g(\la)$ satisfies
(\ref{dg=rho}).

\begin{remark} Formula (\ref{form:xi}) implies that the image of 
$\{-,-\}_G$ under the map $\log : G \to \g$ is given by 
$$
\{f,g\}(x) = \langle df(x) \otimes dg(x),
\big( \on{ad}(x) \otimes \big( {1\over 2} \on{ad}(x) 
\on{coth}({1\over 2}\on{ad}(x))\big)  \big)(t) 
+ \on{ad}(x)^{\otimes 2}(r_0) \rangle ; 
$$
this is a result of \cite{FM2}. 
\end{remark}

\subsection{Constructions based on quantization of Lie bialgebras}
\label{sect:quant}

We now give a proof of Thm. \ref{thm:main} based on 
the theory of quantization of Lie bialgebras (Sect. \ref{section1:2.5}), and then 
give proofs of Props. \ref{prop:1} and \ref{prop:2} based on the same theory
(Sects. \ref{subsect:toto} and \ref{subsect:titi}). 

\subsubsection{Unconditional proof of Thm. \ref{thm:main}}\label{section1:2.5}

It follows from \cite{EH} that the subalgebras $(U^J)'$ and $U'$ of
$U$ are equal. According to 
\cite{Dr:QG,Gav}, $(U^J)'$ is a flat deformation of $\cO_{G^*} := 
U(\g^*)^*$. On the other hand, $U'$ is a flat deformation of 
$\cO_{\g^*} := \wh S(\g)$. So the equality $i_\hbar : (U^J)' \to U'$ 
induces an isomorphism of Poisson algebras 
$i : \cO_{G^*} \stackrel{\sim}{\to} \cO_{\g^*}$, and 
therefore a Poisson isomorphism $\g^* \to G^*$. 

\subsubsection{Construction of $g(\lambda)$ (proof of Prop. \ref{prop:1})}\label{subsect:toto}

Set $U := U(\g)[[\hbar]]$, $U' := \{x\in U | \forall n\geq 0, 
\delta_n(x)\in \hbar^n U^{\wh\otimes n}\} 
= U(\hbar\g[[\hbar]]) \subset U$ (here $\delta_n = (\id - \eta\circ
\eps)^{\otimes n} \circ \Delta^{(n)}$, $\eta$ is the unit map $\kk[[\hbar]]
\to U$). 
As a $\kk[[\hbar]]$-algebra, $U'$ is a flat deformation of 
$\wh S(\g) = \kk[[\g^*]]$. 

Let $\Phi\in U^{\wh\otimes 3}$ be an admissible associator. 
This means that 
$$
\Phi \in 1 + {{\hbar^2}\over {24}}[t^{1,2},t^{2,3}] 
+ \hbar^2 U^{\wh\otimes 3}, 
\quad 
\hbar\log(\Phi) \in (U')^{\wh\otimes 3}, 
$$
$(U,m,\Delta_0,\cR_0 = 1, \Phi)$ is a quasitriangular quasi-Hopf algebra. 
We also require $\eps^{(i)}(\Phi)=1$, $i=1,2,3$ 
(here $\eps^{(1)} = \eps\otimes
\id \otimes \id$, etc., $\eps$ is the counit, $m$ and $\Delta_0$ are the
undeformed product and coproduct of $U$). 
According to \cite{EH}, any universal Lie associator 
gives rise to an admissible associator. 

According to \cite{EK}, there exists a twist killing $\Phi$, and 
according to \cite{EH}, this twist can then be made admissible 
by a suitable gauge transformation. The resulting twist 
$J$ satisfies the following 
conditions: $J\in U^{\wh\otimes 2}$, $J = 1 - \hbar r/2 + o(\hbar)$, 
$\hbar\log(J)\in (U')^{\wh\otimes 2}$, $(\eps\otimes \id)(J) = 
(\id\otimes \eps)(J) = 1$, 
\begin{equation} \label{J:Phi}
\Phi = (J^{2,3} J^{1,23})^{-1} J^{1,2} J^{12,3}. 
\end{equation}
Then $U^J := (U,m,\Delta^J,\cR)$ is a quasitriangular Hopf algebra
quantizing $(\g,r)$. Here $\Delta^J(x) = J\Delta_0(x)J^{-1}$ and 
$\cR = J^{2,1} e^{\hbar t/2} J^{-1}$. 

We have $\on{Ker}(\eps) \cap U' \subset \hbar U$, therefore 
$\log(J)\in U \wh\otimes U'$. 
Then its reduction mod $\hbar$, denoted $g(\lambda) = g^{1,2} 
= \log(J)_{|\hbar=0}$, 
belongs to $U(\g) \wh\otimes \wh S(\g) = U(\g)[[\g^*]]$ 
(formal series on $\g^*$ with coefficients in $U(\g)$). 
The reduction mod $\hbar$ of (\ref{J:Phi}) is 
$g^{12,3} = g^{1,3}g^{2,3}$. Since we also have $(\eps\otimes \id)(g)=1$, 
we get $g = \exp(A)$, with $A\in \g\otimes \wh S(\g)_{>0}$. 

\begin{lemma} \label{lemma:basic}
$g(\lambda)$ satisfies (\ref{dg=rho}). 
\end{lemma}

{\em Proof.} According to \cite{EE}, $\Phi\in U^{\wh\otimes 2}
\wh\otimes U'$ has the expansion $1+\hbar \phi_1 + o(\hbar^2)$, 
where $\phi_1 \in U^{\wh\otimes 2} \wh\otimes U'$ is such that 
$(\phi_1 - \phi_1^{2,1,3})_{|\hbar=0} = -\rho_{\on{AM}}$. 

If $x\in U(\g)\wh\otimes \wh S(\g)$, we denote by $\bar x$ a lift of $x$ in 
$U\wh\otimes U'$. 

Let us expand  $\log(J)$ as $\bar A + \hbar A_1 + o(\hbar)$, with $A_1\in 
U\wh\otimes U'$. 

Then $J^{1,23} = \exp(\bar A^{1,3} + \hbar(A_1^{1,3} + d_2 A^{1,3}) 
+ o(\hbar)) = J^{1,3} (1 + \hbar \overline{g_1^{-1}d_2(g_1)} + o(\hbar))$. 

We have $[J^{1,3},J^{2,3}] = \hbar \overline{\{g^{1,3},g^{2,3}\}} + o(\hbar)$, 
so $(J^{12,3})^{-1}[J^{1,3},J^{2,3}] = \hbar \overline { 
(g^{1,3}g^{2,3})^{-1}\{g^{1,3},g^{2,3}\} } + o(\hbar) = 
\hbar \overline{ \langle \id\otimes \id\otimes \lambda, 
[g_1^{-1}d_3(g_1),g_2^{-1}d_3(g_2)] \rangle} + o(\hbar)$. 
So we get 
$$
J^{12,3} = J^{2,3} J^{1,3} (1 + \hbar\psi_1 + o(\hbar)), 
$$
where $\psi_1 \in U^{\wh\otimes 2}\wh\otimes U'$ is such that 
$(\psi_1 - \psi^{2,1,3})_{|\hbar=0} = \langle \id\otimes \id\otimes \lambda,
[g_1^{-1}d_3(g_1),g_2^{-1}d_3(g_2)] \rangle$. 

Then (\ref{J:Phi}) gives 
$$
1 + \hbar \phi_1 + o(\hbar) = (1 - \hbar \overline{g_1^{-1} d_2(g_1)} 
+ o(\hbar))
(J^{1,3})^{-1} (J^{2,3})^{-1} (1 - \hbar r/2 + o(\hbar)) J^{2,3}
J^{1,3} (1 + \hbar\psi_1 + o(\hbar)). 
$$
The reduction modulo $\hbar$ of $(J^{1,3})^{-1} (J^{2,3})^{-1} r 
J^{2,3} J^{1,3}$ is $\Ad(g\otimes g)^{-1}(r)$. Then substracting $1$, 
dividing by $\hbar$, reducing modulo $\hbar$ and 
antisymmetrizing the two first tensor factors, we get the lemma. 
\hfill \qed \medskip

More generally, assume that $J'\in U^{\wh\otimes 2}$ satisfies 
$J' = 1 - \hbar r/2 + o(\hbar)$, $\hbar \log(J') \in (U')^{\wh\otimes 2}$, 
$(\eps\otimes \id)(J')=1$ and 
\begin{equation} \label{eq:J'}
\Phi = (J^{\prime 2,3} J^{\prime 1,23})^{-1} J^{1,2} J^{\prime 12,3}. 
\end{equation}
Then $J'\in U\wh\otimes U'$, and its reduction $g'$ modulo $\hbar$
satisfies: $g'\in \exp(\g\otimes \wh S(\g)_{>0})$, and equation (\ref{dg=rho})
with $g$ replaced by $g'$. Moreover, there exists $u\in U^\times$ such that 
$J' = u^2 J (u^{12})^{-1}$. 

\begin{remark} Equation (\ref{eq:J'}) can be interpreted as saying that 
$J'$ is a vertex-IRF transformation relating $J^{1,2}$ 
and $\Phi$, and equation (\ref{dg=rho}) for $g'$ is the classical limit of 
this statement (see \cite{EN}). Vertex-IRF transformations are a special 
kind of non-invariant dynamical gauge transformations, which map a constant, 
but non-invariant twist to an invariant, but non-constant 
(i.e., dynamical) one.  
\hfill \qed \medskip 
\end{remark}

Let $U'_0 := \Ker(\eps)\cap U'$. Then $\hbar^{-1}U'_0 \subset U$
is a Lie subalgebra for the commutator. This Lie algebra acts on the set of 
solutions of (\ref{eq:J'}) by $\delta_u(J) = u^2 J - J u^{12}$ ($\delta_u(J)
\in U\wh\otimes U'$ because it is equal to $[u^2,J] - J(u^{12}-u^2)$); this 
means that if $\eps$ is a formal parameter with $\eps^2=0$, then $(\id + \eps
\delta_u)(J')$ is a solution of (\ref{eq:J'}) if $J'$ is. 

The reduction modulo $\hbar$ of this action may be described as follows. 
The Lie algebra $(\wh S(\g)_{>0},\{-,-\})$ acts on the set of solutions of 
(\ref{dg=rho}) by $\delta_f(g) = \{1\otimes f,g\} - g \cdot df$, i.e., action
(\ref{action}). When restricted to the Lie subalgebra $\wh S(\g)_{>1}$, 
this action is the infinitesimal of the right action of 
$\on{Map}^{ham}_0(\g^*,G)$ on the set of solutions of (\ref{dg=rho}), 
given by $(g*\alpha)(\lambda) = g(\Ad^*(\alpha(\lambda))(\lambda)) 
\alpha(\lambda)$. 

\subsubsection{Isomorphism $\g^* \to G^*$ (proof of Prop. \ref{prop:2})} 
\label{subsect:titi}

Let us construct the maps  $a : \cO_G \to \cO_{\g^*}$ and $b : \cO_G \to \cO_{G^*}$
out of the quantization of Lie bialgebras. 

We first construct quantized versions of $a$ and $b$. 
We have $\hbar\log(J)\in (U')^{\otimes 2}$, therefore 
$J,J^{-1}\in U\wh\otimes U'$ and $Je^{\hbar t}J^{-1} \in U\wh\otimes U'$. 
We then define $a_\hbar : U^* \to U'$ by 
$a_\hbar(f) := (f\otimes \id)(Je^{\hbar t}J^{-1})$. 
Now we also have $J^{2,1}\in U\wh\otimes U'$, hence 
$\cR,\cR^{2,1}\in U\wh\otimes U' = U\wh\otimes (U^J)'$. 
We then define the linear map $b_\hbar : U^* \to (U^J)'$ 
by $b_\hbar(f) := (f\otimes \id)(\cR^{2,1}\cR)$. 
Then $a_\hbar \circ b_\hbar^{-1} : U' \to (U^J)'$ coincides with 
$i_\hbar$. We define $a,b$ as the reductions modulo $\hbar$ of 
$a_\hbar,b_\hbar$. 

Let us now compute the classical limit of $a_\hbar$. Define 
maps $j_\hbar,\iota_\hbar,j'_\hbar : U^* \to U'$ by $j_\hbar(f) = 
(f\otimes \id)(J)$, $j'_\hbar(f) = (f\otimes \id)(J^{-1})$. 
Recall that the reductions modulo $\hbar$ of all three elements 
$J,J^{-1},e^{\hbar t}$ in $U(\g) \wh\otimes \wh S(\g)$ are of the form $K =
\exp(k)$, where $k\in \g\otimes \wh S(\g)_{>0}$.  

\begin{lemma} If $K\in U(\g) \wh\otimes \wh S(\g)$ is of the form 
$\exp(k)$, where $k\in \g\otimes \wh S(\g)_{>0}$, then the 
morphism $\cO_G = U(\g)^* \to \wh S(\g) = \cO_{\g^*}$ given by $f\mapsto 
(f\otimes \id)(K)$ is dual to the morphism $\g^* \to G$, 
$\lambda \mapsto e^{k(\lambda)}$. 
\end{lemma}

{\em Proof.} We compose this morphism with the 
transpose of the inverse of the symmetrization map 
${}^t\on{Sym}^{-1}: S(\g)^* \to U(\g)^*$. The morphism 
${}^t\on{Sym}^{-1}$ corresponds to the logarithm map $G\to \g$. 
Now the composed morphism $S(\g)^* \to \wh S(\g)$ is given by 
$f\mapsto (f\otimes \id)((\on{Sym}^{-1} \otimes \id)(K))$. Now 
$(\on{Sym}^{-1}\otimes \id)(K)$ is $\exp(k)$, where the 
exponential is now taken 
in $S(\g)\wh\otimes \wh S(\g)$. The morphism $S(\g)^* \to \wh S(\g)$, 
$f\mapsto (f\otimes \id)(\exp(k))$, is an algebra morphism, 
taking the function $X\mapsto \alpha(X)$ on $\g$ ($\alpha\in\g^*$) 
to the function $\lambda\mapsto \alpha(k(\lambda))$ on $\g^*$, 
and therefore corresponds to the morphism $\g^*\to\g$, 
$\lambda\mapsto k(\lambda)$. Composing with the exponential, we get the
announced morphism.     
\hfill \qed \medskip 

It follows that the reductions modulo $\hbar$ of the morphisms 
$j_\hbar,\iota_\hbar,j'_\hbar$ are morphisms $\cO_G \to \cO_{\g^*}$, 
corresponding to morphisms $j,\iota,j' : \g^* \to G$ such that 
$j(\la) = g(\la)$, $\iota(\la) = e^{\la^\vee}$, $j'(\la)
= g(e^\la)^{-1}$.  

Then $a_\hbar = m^{(2)} \circ (j_\hbar\otimes \iota_\hbar\otimes 
j'_\hbar)\circ \Delta^{(2)}$, so $a : \cO_{G^*} \to \cO_G$ corresponds 
to the composed map $G^* \stackrel{\on{diag^{(2)}}}{\to} (G^*)^3 
\stackrel{(j,\iota,j')}{\to} G^3 \stackrel{\on{product}^{(2)}}{\to} G$. 
Therefore we get $a(\la) = g(\la) e^{\la^\vee} g(\la)^{-1}$.  

We now compute the classical limit of $b_\hbar$.  
Define maps $L_\hbar, R'_\hbar : (U^J)^* \to (U^J)'$ by 
$L_\hbar(f) := (f\otimes \id)(\cR^{2,1})$, 
$R'_\hbar(f) := (f\otimes \id)(\cR)$. Then $L_\hbar$ is an 
antimorphism of algebras and a morphism of coalgebras, while 
$R'_\hbar$ is a morphism of algebras and antimorphism of 
coalgebras. Their reductions $L,R'$ modulo $\hbar$ are morphisms 
$\cO_{G^*} \to \cO_G$ (anti-Poisson coalgebra morphism in the case of $L$, 
Poisson anti-coalgebra in the case of $R'$). Using \cite{EGH}, 
appendix, one shows that these morphisms correspond to 
the morphisms of formal groups $L,R' : G^* \to G$ (antimorphism 
in the case of $R'$), corresponding to the morphisms 
$L,R' : \g^* \to \g$, given by $L(\lambda) := (\la\otimes\id)(r)$ and 
$R'(\la) := (\la\otimes \id)(r^{2,1})$. Here $L$ is a Lie algebra, 
anti-Lie coalgebra morphism and $R'$ is an anti-Lie algebra, 
Lie coalgebra morphism. Now $b_\hbar = m \circ (L_\hbar \otimes 
R_\hbar) \circ \Delta$, so $b : \cO_{G^*} \to \cO_G$ corresponds to 
the composed map $G^* \stackrel{\on{diag}}{\to} (G^*)^2 
\stackrel{(L,R')}{\to} G^2 \stackrel{\on{product}}{\to} G$, i.e., 
$g^* \mapsto L(g^*)R'(g^*) = L(g^*)R(g^*)^{-1}$. 

Finally, the isomorphism $i : \g^* \to G^*$ is equal to 
$b^{-1} \circ a$, so it takes $\lambda\in \g^*$ to $g^* \in G^*$
such that 
$$
L(g^*)R(g^*)^{-1} = g(\lambda) e^{\lambda^\vee} g(\lambda)^{-1}. 
$$
So it coincides with the isomorphism obtained in Prop. \ref{prop:2}. 

\subsection{On the group $\on{Map}^{ham}_0(\g^*,G)$}

In this section, we give a "quantum" proof of the following statement, 
which was used in Sect. \ref{sect:geom} (and can also be proved in the 
setup of this section).  

\begin{proposition} $\on{Map}^{ham}_0(\g^*,G)$ is a subgroup of 
$\on{Map}_0(\g^*,G)$. 
\end{proposition}

{\em Proof.} Set $U'_0 = \on{Ker}(\eps) \cap U'$. 
The map of $U'\to \wh S(\g)$ of reduction by $\hbar$
restricts to $U'_0\to \wh S(\g)_{>0}$, so we get a map 
$U'_0 \to \g$. 

Consider the set of all $J'$, such that 
$\hbar\on{log}(J') \in (U'_0)^{\wh\otimes 2}$, 
the image of $\hbar\on{log}(J')$ in $\g^{\otimes 2}$
is zero, and $J^{\prime 2,3} J^{\prime 1,23} = J^{\prime 1,23}$. 

This is the set of elements of the form $J'(u) = u^2 (u^{12})^{-1}$, 
where $u \in \on{exp}(U'_0/\hbar)$ is such that the image of 
$\on{log}(u)$ under $\wh S(\g)_{>0} \to \g$ is zero. 
Indeed, one recovers $u$ from $J'$ by 
$u = (\id\otimes\eps)(J')^{-1}$. 

We denote by $A(\la)$ the image of 
$\hbar\on{log}(J')$ by $\wh S(\g)_{>0} \wh\otimes 
\wh S(\g)_{\geq 2} \to \g\otimes \wh S(\g)_{\geq 2}$, and set $g(\la)
:= \exp(A(\la))$. Then $g(\la)\in \on{Map}_0(\g^*,G)$.  

As in Lemma \ref{lemma:basic}, one proves that $g(\la)$ satisfies equation 
(\ref{eq:map:ham}). Hence we get a map of sets $\{J'$ as above$\}
\to \on{Map}^{ham}_0(\g^*,G)$, which is surjective. 

On the other hand, the set of all $J'(u)$ is equipped with a product, 
such that $J'(u) * J'(v) = J'(uv)$. This product expresses as follows: 
$J'(u) * J'(v) = \big( u^2 J'(v) (u^2)^{-1}\big) J'(u)$; the classical 
limit of this expression is the product formula for $\on{Map}_0(\g^*,G)$. 
So $\on{Map}^{ham}_0(\g^*,G)$ is a subgroup of $\on{Map}_0(\g^*,G)$.
\hfill \qed \medskip

\section{Quantization of $\rho_{\on{FM}}$} 

In this section, we prove Thm. \ref{thm:2}. 

Let $\Phi_{univ}$ be a universal Lie associator defined over $\kk$
with parameter $\mu=1$. So $\Phi_{univ} = 1 + {1\over {24}}
[t_{12},t_{23}] +$ terms of degree $>2$, where $t_{ij}$ is the universal
version of $t^{i,j}$ (see \cite{Dr:Gal}). Set 
$\Phi_\nu := \Phi_{univ}(2\hbar\nu t^{1,2},
2\hbar\nu t^{2,3})$ and $\Phi := \Phi_{1/2}$. Then 
$$
(U(\g)[[\hbar]], m, \Delta_0, \Phi_\nu)
$$  
is a quasi-Hopf algebra; its classical limit is the quasi-Lie bialgebra 
$(\g,\delta = 0,\nu^2[t^{1,2},t^{2,3}])$. Let $J$ be an admissible 
twist killing $\Phi$, and let us twist this quasi-Hopf algebra by $J$. 
We obtain the quasi-Hopf algebra 
$$
(U(\g)[[\hbar]], m, \Delta_0^J, \Phi_\nu^J), 
$$  
where $\Delta_0^J(x) = J \Delta_0(x) J^{-1}$ and 
$\Phi_\nu^J = J^{2,3}J^{1,23}
\Phi_\nu (J^{-1})^{12,3} (J^{-1})^{1,2}$. Its classical limit is 
the quasi-Lie bialgebra
$(\g,\delta(x) = [x^1 + x^2,r], Z_\nu)$. 

Now $(U(\g)[[\hbar]], m, \Delta_0^J)$ is a Hopf algebra 
quantizing $(\g,\delta)$, which we denote by $U_\hbar(\g)$, 
and $\Phi_\nu^J$ is an associator for this quantized universal enveloping
algebra, with classical limit $Z_\nu$. 

$\Phi_\nu^J$ clearly satisfies the invariance and pentagon equations.
We have $J = 1 + \hbar j_1$, where $j_1\in U\wh\otimes U'$, so 
$J^{2,3},J^{1,23},(J^{-1})^{1,23}$ and $(J^{-1})^{1,2}$ are all of the form 
$1 + \hbar k$, $k\in U^{\wh\otimes 2} \wh\otimes U'$. Hence $\Phi_\nu^J = 1 +
\hbar \psi_\nu$, where $\psi_\nu\in U^{\wh\otimes 2} \wh\otimes U'$. 
Set $U_\hbar := U_\hbar(\g)$, then $U_\hbar' = U'$, so 
$\psi_\nu\in U_\hbar^{\wh\otimes 2} \wh\otimes U_\hbar'$. 
We now compute $(\psi_\nu - \psi_\nu^{2,1,3})_{|\hbar=0}$; 
this is an element of $U(\g)^{\otimes 2} \wh\otimes \cO_{G^*}
\simeq U(\g)^{\otimes 2} \wh\otimes \cO_{\g^*}$.  

We have $\Phi_\nu^J = \Phi^J + (\Phi_\nu - \Phi)^J = 1 + 
\hbar(\Phi_\nu - \Phi)^J$. Let us define $\phi,\phi_\nu$ by 
$\Phi = 1 + \hbar \phi$, $\Phi_\nu = 1 + \hbar \phi_\nu$, then 
$(\psi_\nu)_{|\hbar=0} = \big( (\phi_\nu - \phi)^J\big)_{\hbar=0}$, 
therefore 
\begin{eqnarray*}
(\psi_\nu - \psi_\nu^{2,1,3})_{|\hbar=0} & 
= 
((\phi_\nu - \phi_\nu^{2,1,3})^J)_{|\hbar=0}
- ((\phi - \phi^{2,1,3})^J)_{|\hbar=0}
\\ & = \Ad(g(\lambda))^{\otimes 2} (\rho_{\on{AM}}(\lambda) - 
\rho_{\on{AM}}^\nu(\lambda)) 
\end{eqnarray*} 
as a formal function $\g^* \to \wedge^2(\g)$. Here 
$\rho_{\on{AM}}^\nu(\lambda) = 2\nu\rho_{\on{AM}}(2\nu\lambda)$. 
Since $\rho_{\on{AM}}$ and $\rho_{\on{AM}}^\nu$ are $G$-equivariant, 
this is 
$$
(\rho_{\on{AM}} - \rho_{\on{AM}}^\nu)(\Ad^*(g(\lambda))(\lambda))
= \Big( \id\otimes \big( {1\over 2}
{{e^{\on{ad}(\bar\lambda^\vee)} + \id}\over
{e^{\on{ad}(\bar\lambda^\vee)} - \id}}
- \nu {{e^{2\nu\on{ad}(\bar\lambda^\vee)} + \id}\over
{e^{2\nu\on{ad}(\bar\lambda^\vee)} - \id}}
\big)\Big)(t)
= \rho_{\on{FM}}(g^*). 
$$ 
Here we set $\bar\lambda := \Ad^*(g(\lambda))(\lambda)$ and 
we use the relation $L(g^*)G(g^*)^{-1} = e^{\bar\lambda^\vee}$. 

We have proved: 

\begin{theorem} \label{thm:end}
The quantized universal enveloping algebra 
$U_\hbar(\g) = U(\g)[[\hbar]]^J$, together with the pair 
$(\bar J,\bar \Phi)$ defined by $\bar J = 
\bar \Phi = \Phi_\nu^J$, 
is a quantization of the Poisson-Lie dynamical $r$-matrix 
$\rho_{\on{FM}}(g^*)$ over $(G^*,\g,Z_\nu)$. 
\end{theorem}

\subsection*{Acknowledgements} We are grateful to A. Alekseev and E. Meinrenken, and to V. Toledano Laredo, for independently pointing out 
the necessity of the nondegeneracy condition of $t$ in the formulation of Prop. \ref{prop:2}, 
which was overlooked in a previous version of this work.


\begin{thebibliography}{EnGH}

\bibitem[Al]{Al}
A. Alekseev, {\it On Poisson actions of compact Lie groups 
on symplectic manifolds,} J. Differential Geometry
{\bf 45} (1997), 241-56.

\bibitem[AM1]{AM}
A. Alekseev, E. Meinrenken, {\it The non-commutative 
Weil algebra,} Invent. Math. {\bf 139} (2000), 135-72.

\bibitem[AM2]{AM2}
A. Alekseev and E. Meinrenken, {\it Linearization of Poisson-Lie structures,} J. Sympl. Geom. {\bf 14} (2016), 227-267, arXiv:1312.1223. 

\bibitem[BDF]{BDF} J. Balog, L. Dabrowski,and L. Feh\'er,
{\it Classical $r$-matrix and exchange algebra in WZNW and 
Toda field theories,}
Phys. Lett. B, {\bf 244} (1990), issue 2, 227-34.

\bibitem[BFP]{BFP}
J. Balog, L. Feh\'er, L. Palla, {\it Chiral extensions of the WZNW phase space,
Poisson-Lie symmetries and groupoids,} preprint hep-th/9910046, Nucl. Phys. B 
{\bf 568} (2000), 503-42. 

\bibitem[Bo]{Bo}
P. Boalch, {\it Stokes matrices, Poisson-Lie groups and Frobenius manifolds,}
Invent. Math.  {\bf 146} (2001), no. 3, 479-506. 

\bibitem[Ch]{Ch}
V. Chloup-Arnould, {\it Linearization of some Poisson-Lie tensor,} 
J. Geom. Phys. {\bf 24} (1997), no. 1, 46-52.

\bibitem[Dr1]{Dr:QG} V. Drinfeld, {\it Quantum groups,} 
Proceedings of the ICM-86 (Berkeley), 
798-820, Amer. Math. Soc., Providence, RI, 1987. 

\bibitem[Dr2]{Dr:QH} V. Drinfeld, {\it Quasi-Hopf algebras,} 
Leningrad Math. J. {\bf 1} (1990), no. 6, 1419-57. 

\bibitem[Dr3]{Dr:Gal} V. Drinfeld, {\it On quasitriangular quasi-Hopf 
algebras and on a group that is closely connected with 
${\rm Gal}(\overline\QQ/\QQ)$,}
Leningrad Math. J. {\bf 2} (1991), no. 4, 829-60. 

\bibitem[EE]{EE} B. Enriquez, P. Etingof, {\it Quantization of 
Alekseev-Meinrenken dynamical $r$-matrices,} preprint math.QA/0302067
(in memory of F.I. Karpelevich), AMS Transl. {\bf 210} (2003), 
no. 2, 81-98. 

\bibitem[EEM]{EEM} B. Enriquez, P. Etingof, I. Marshall, {\it Quantization of 
some Poisson-Lie dynamical $r$-matrices and Poisson homogeneous spaces,} 
preprint math.QA/0403283 (in memory of J. Donin).   

\bibitem[EGH]{EGH} B. Enriquez, F. Gavarini, G. Halbout, {\it Uniqueness of
braidings of quasitriangular Lie bialgebras and lifts of classical 
$r$-matrices,} Internat. Math. Res. Notices, {\bf 46} (2003), 2461-86. 

\bibitem[EH]{EH} B. Enriquez, G. Halbout, {\it Poisson algebras 
associated to quasi-Hopf algebras,} Adv. Math.  {\bf 186}  (2004),  no. 2, 
363-95. 

\bibitem[EK]{EK} P. Etingof, D. Kazhdan, {\it Quantization of Lie 
bialgebras, I, II,} Selecta Math. (N.S.) \textbf{2} (1996), no. 1, 1-41;
\textbf{4} (1998), no. 2, 213-31. 

\bibitem[EN]{EN} P. Etingof, D. Nikshych, 
{\it Vertex-IRF transformations and quantization of 
dynamical $r$-matrices,} Math. Res. Lett.  {\bf 8}  (2001),  
no. 3, 331-45.

\bibitem[ES]{ES} P. Etingof, O. Schiffmann, 
{\it On the moduli space of classical dynamical $r$-matrices,}
Math. Res. Lett. {\bf 8} (2001), no. 1-2, 157-70.

\bibitem[FM1]{FM1} L. Feh\'er, I. Marshall, {\it On a Poisson-Lie analogue of 
the classical dynamical Yang-Baxter equation for self-dual Lie algebras,} 
Lett. Math. Phys. {\bf 62} (2002), 51-62. 

\bibitem[FM2]{FM2} L. Feh\'er, I. Marshall, {\it The non-abelian momentum 
map for Poisson-Lie symmetries on the chiral WZNW phase space,} 
Int. Math. Res. Not. {\bf 49} (2004), 2611-36. 

\bibitem[Gav]{Gav} F. Gavarini, {\it The quantum duality principle,}
Ann. Inst. Fourier (Grenoble) {\bf 52} (2002), no. 3, 809-34.

\bibitem[GW]{GW} 
Viktor L. Ginzburg, A. Weinstein, {\it Lie-Poisson structure on some 
Poisson-Lie  groups,} J. Amer. Math. Soc.  5:2  (1992), 445-53.

\bibitem[STS]{STS} M. Semenov-Tian-Shansky, 
{\it Dressing transformations and Poisson group actions,} 
Publ. Res. Inst. Math. Sci. {\bf 21} (1985), no. 6, 1237-1260.



\end{thebibliography}
\end{document}